\documentclass{amsart}

\newtheorem{theorem}{Theorem}

\newtheorem{corollary}{Corollary}
\theoremstyle{remark}
\newtheorem{remark}{Remark} 

\newcommand{\aut}[1]{\operatorname{Aut}(#1)}
\newcommand{\rank}{\operatorname{rank}}
\newcommand{\thmfont}{\sl}
\newcommand{\integer}{{\mathbb Z}}
\newcommand{\complex}{{\mathbb C}}
\newcommand{\cstar}[1]{(\complex^*)^{#1}}

\newcommand{\ox}{o_{\text{\rm excl}}}
\newcommand{\Ox}{O_{\text{\rm excl}}}
\def\spa#1{\complex^#1}
\newcommand{\cal}{\mathcal}

\newcommand{\p}[2]{\frac{\partial #1}{\partial #2}}
\renewcommand{\a}{\alpha}
\newcommand{\z}{\zeta}
\catcode`"=\active
\def"#1"{\ref{#1}}

\begin{document}

\title{Complete vector fields on $\cstar{n}$}

\author{Erik Anders\'en}


\address{Department of Mathematics\\ Purdue University\\1395
Mathematical Sciences Bldg.\\West Lafayette, IN 47907-1395\\USA}

\email{andersen@@math.purdue.edu}

\thanks{Supported by a grant from STINT}

\begin{abstract}
We prove necessary and sufficient conditions for a rational
vector field on $(\complex^*)^n$ to be complete.
\end{abstract}


\maketitle

\section{Introduction} This paper has been motivated by the desire to
understand the group of holomorphic automorphisms of certain complex
manifolds $M$. In a number of cases the group, $\aut M$, is known to
be a finite dimensional Lie group: this is notably so when $M$
is a bounded domain in $\spa{n}$ or, more generally, a hyperbolic
manifold; or when $M$ is compact, see 
\cite{kobayashi:hyperbolic}. The manifolds that we are concerned with are not
such. Little is known about the automorphism groups of nonhyperbolic
affine varieties beyond the fact that they can be huge
(\cite{varolin:density}); but for the
case $M=\spa{n}$ see
\cite{andersen:volume-preserving,andersen-lempert:automorphisms,rosay-rudin:maps,forstneric:actions}.

If $2n$ generic planes are removed from
$\spa{n}$ the resulting space is known to be hyperbolic (see
\cite{bloch:systemes,green:hyperbolicity}), and therefore the
automorphism groups are finite dimensional.
The case of $\spa{n}$ minus $m$ hyperplanes
($n<m<2n$) is subject to research by the author.
In this
paper we take $M$ to be $\spa{n}$ minus $n$ hyperplanes in general
position ($n\ge2$). If $\complex^*=\complex\setminus\{0\}$, we can take
$M$ to be $\cstar{n}$. These manifolds are known to have infinite
dimensional automorphism groups. In the case $n=2$ we can get
automorphisms of $\cstar{2}$ by taking an arbitrary holomorphic
function $f$, two 
integers $n_1,n_2$ and  complex numbers $c_1,c_2$ and forming
 \begin{equation}\label{c-2-automorphism}
(z_1,z_2)\mapsto
(z_1e^{n_2f(z_1^{n_1}z_2^{n_2})},z_2e^{-n_1f(z_1^{n_1}z_2^{n_2})}).
 \end{equation}
Once we observe that this mapping preserves $z_1^{n_1}z_2^{n_2}$ we
easily see that it is bijective. This verifies the claim that the
automorphism groups of $\cstar{n}$ are infinite dimensional. Other
automorphisms of $\cstar{2}$ are given by 
\begin{equation}\label{matrix}
(z_1,z_2)\mapsto(z_1^{a_{11}}z_2^{a_{12}},z_1^{a_{21}}z_2^{a_{22}}),
\end{equation}
where the integers $a_{ij}$ satisfy $a_{11}a_{22}-a_{12}a_{21}=1$.
It is conjectured that these mappings generate the full automorphism
group of $\cstar2$. Nishimura \cite{nishimura:automorphisms} proves
that any automorphism of $\cstar{2}$ which extends to $\complex^2$ and
preserves the volume form $dz_1\wedge dz_2$ is of the form
("c-2-automorphism"), with $n_1=n_2=1$. He also has results about
automorphisms of 
$\complex\times\complex^*$ --- see \cite{nishimura:applications}.

Peschl \cite{peschl:automorphismes} claimed to have proved that all
automorphisms of 
$\cstar2$ that extend to mappings of $\spa2$ preserve the volume form
 \[
\frac{dz_1\wedge dz_2}{z_1z_2},
 \]
but his proof has been found to be incomplete. Accordingly, there is a
conjecture that any automorphism of $\cstar2$ preserves this form.

An action of $\complex$ on $M$ is a family $\{\phi_t:t\in\complex\}$ of
automorphisms of $M$ such that $\phi_s\circ\phi_t=\phi_{s+t}$.
Suzuki \cite{suzuki:operations,suzuki:groupe-additif} has studied
actions on
a two dimensional manifold. He proves among other things that if an
action of $\complex$ on $\cstar{2}$ factors through an action of
$\complex^*$ on $\cstar{2}$ then the action is by linear mappings
$(z_1,z_2)\mapsto(e^{n_2t}z_1,e^{-n_1t}z_2)$. 

An
action is generated by a complete holomorphic vector field. (A
holomorphic vector field is said to be complete if all integral curves
are entire holomorphic functions, see \cite{forstneric:actions}.)
In this paper we
give a complete characterization of complete holomorphic rational vector 
fields on $\cstar{n}$. Such fields can be written
 \[
\dot z_i=z_ip_i(z)
 \]
where $z=(z_1,\dots,z_n)$ and $p_i$ are Laurent polynomials in
$z_1,\dots,z_n$,  that is, polynomials in these variables and their
inverses.  As a corollary, we
prove the above conjecture for automorphisms coming from 
such fields. 

\section{Preliminaries}
We will need the following elementary facts from Nevanlinna theory
(see e.g. \cite{lang:hyperbolic}). For any
meromorphic function of one complex variable $f(z)$ we set
\begin{gather*}
m(r,f)=\frac1{2\pi r}\int_{|z|=r}\log^+|f(z)|\,|dz|,\\
N(r,f)=\sum_{f(a)=\infty,\ a\ne0}\log^+\left|\frac ra\right|+k\log r
\quad\text{and}\\
T(r,f)=m(r,f)+N(r,f),
\end{gather*}
where $k$ is the order of the pole of $f$ at $0$. ($k=0$ if $f$ is
regular at $0$.)
In the definition of $N$, the sum is taken over all poles $a$ of $f$, with
regard to multiplicity. The function $T(r,f)$ is called the Nevanlinna 
characteristic of $f$. The following properties
are easily verified.
\begin{equation} \label{Nevanlinna-rules}
\begin{gathered} 
T(r,f_1+f_2)\le T(r,f_1)+T(r,f_2)+O(1)\\
T(r,f_1f_2)\le T(r,f_1)+T(r,f_2),\\
T(r,f^d)=dT(r,f)\quad(d>0).
\end{gathered}
\end{equation}
Here and in the sequel, the estimates $O(g)$ and $o(g)$ are as $r\to\infty$.
The first fundamental theorem of Nevanlinna theory says that
\[
T(r,1/f)=T(r,f)+O(1).
\]
We also have the Lemma of the logarithmic derivative (LLD),
\[
m(r,f'/f)=\ox(T(r,f)).
\]
Here $\ox$ means that the estimate holds outside a set
of finite measure. We also use the corresponding notation $\Ox$. From LLD
and the preceding inequalities it follows that
\begin{gather}
m(r,f^{(k}))\le(1+\ox(1)))T(r,f)\\
m(r,f^{(k)}/f)=\ox(T(r,f))
\end{gather}
for all positive integers $k$.

\section{Borel's Theorem} We need a version of a classical theorem of
Borel. Since it is generally stated in a slightly different form (the
function $f$ below is usually $0$ or $1$, and the conclusions are also
slightly different from what we need) I include a complete proof. For
Borel's original theorem, see \cite{borel:fonctions}.

\begin{theorem}\thmfont Let $f$ and $u_1,\dots,u_n$ be entire
functions of one variable satisfying
\begin{equation}\label{u-sum}
\sum_i u_i=f.
\end{equation}
If $u_1,\dots,u_n$ have no zeros then one of the following cases holds.
\begin{enumerate}
\item $T(r,u_i)=\Ox(T(r,f)+1)$ for all $i$.
\item Some non-empty subsum $\sum_{i\in I}u_{i}=0$.
\end{enumerate}
\end{theorem}

\begin{proof} The proof is by induction $n$. If $n=1$ then Case~1 holds
automatically  so there is nothing to prove. We assume that the
theorem holds for sums with less than $n$ terms.

If two terms in ("u-sum") are proportional, say $u_1=au_2$, we can
lump them together  to get a shorter sum
\begin{equation}
(a+1)u_2+u_3+\dots+u_n=f.
\end{equation}
If $a=-1$ then case~(2) holds. Otherwise we apply the theorem to this
shorter sum and whether we get conclusion~(1) or~(2) for this sum we
get the same conclusion for the original sum ("u-sum"). We therefore
assume  that there are no proportional terms in the sum.

If we differentiate (\ref{u-sum})  we get
\[
\sum_i \frac{u_i^{(j)}}{u_i}u_i=f^{(j)}.\quad\text{for $j=0,\dots,n-1$}
\]
This is a linear system for $u_i$. Two cases are possible.
\begin{enumerate}
\item $\det(u_i^{(j)}/u_i)\not\equiv0$. Then we can solve for $u_i$ 
and get
\begin{multline*}
T(r,u_i)=O(T(r,f))+\sum_j O(T(r,u_i^{(j)}/u_i))+O(1)\\
=O(T(r,f))+\sum_i\ox(T(r,u_i))+O(1)
\end{multline*}
by LLD. We therefore have Case~1.
\item $\det(u_i^{(j)}/u_i)\equiv0$. This means that the Wronskian
$W(\{u_i\})=\det(u_i^{(j)})\equiv0$.
The theory of ordinary differential equations now says there are constants
$c_i$ such that
\begin{equation}\label{u dependence}
\sum_i c_iu_i=0.
\end{equation}
We choose the shortest possible such sum, that is, the sum with the smallest
number of non-zero $c_i$.
Possibly after a reordering and a scaling
we may assume that $c_1=-1$, so that $u_1=\sum_{i>1} c_iu_i$. Division by
$u_1$ gives
\begin{equation}\label{v-sum}
\sum_{i>1}c_iu_i/u_1=1.
\end{equation}
We set $v_i=c_iu_i/u_1$ and apply the theorem to the expression (\ref{v-sum}),
which has less than $n$ terms. There are two cases.
\begin{enumerate}
\item $T(v_i,r)=\Ox(T(r,1)+1)=\Ox(1)$ for all $i$. Then all $v_i$ are constant
and $c_iu_i=a_iu_1$ for some $a_i$. We assumed that there were no
proportional terms so this case is excluded.
\item $\sum_{i\in I}v_i=0$ for some set $I$. Then $\sum_{i\in I}c_iu_i=0$ and
this sum is shorter than (\ref{u dependence}). This is a contradiction. \qed

\end{enumerate}

\end{enumerate}
\renewcommand{\qed}{}
\end{proof}

\section{Vector fields}

We start with the notation.  Let $p_i$ be Laurent polynomials and
 \begin{equation}\label{ode}
\dot z_i=z_ip_i(z)\quad i=1,\dots,n
 \end{equation}
be a vector field on $(\complex^*)^n$. Write $p_i(Z)=\sum_\alpha
 p_{i,\a}Z^\a$ for each $i$. Let 
 ${\cal M}$ be the multiplicative group generated by 
$\{Z^\a:p_{i,\a}\ne0\text{ for some }i\}$. $\cal M$ is isomorphic to a
 lattice in $\integer^n$ under the mapping $\integer^n\ni\alpha\mapsto
 Z^\alpha\in\cal M$. As such it has rank at most~$n$. Let $\rank\cal M=m$ and
 $W_i=\prod Z_j^{a_{ij}}$ for 
$i=1,\dots,m$ be a basis for $\cal M$. There are Laurent-polynomials $f_i$ such
that $p_i(Z)=f_i(W)$, where $W=(W_1,\dots,W_m)$. The following theorem
tells us when the field is complete.

\begin{theorem}\label{ode-theorem}\thmfont
Notation as above, we have two cases.

\begin{enumerate}
\item $m=n$. Then (\ref{ode}) is not complete.
\item $m<n$. Then ("ode") is complete if and only if
 \begin{equation}\label{new-ode}
\dot w_i=w_i\sum_j a_{ij}f_j(w), \quad i=1,\dots,m
 \end{equation}
is complete.
\end{enumerate}

\end{theorem}

\begin{proof} We first deal with the case $m=n$. We assume that
(\ref{ode}) is complete and want to derive a contradiction.
Choose $c$ so that $\sum_{\a\in A} p_{i,\a} c^\a\ne0$ for
all non-zero subpolynomials and all indices~$i$. Let $z(t)$ be the
integral curve with $z(0)=c$. Apply Borel's theorem to $\sum_\a p_{i,\a}z^\a=
\dot z_i/z_i$. Because of our choice of initial condition, Case~2 does
not hold. We  therefore have Case~1 and
 \[
T(r,z^\a)=\Ox(T(r,\dot z_i/z_i)+1)=\ox(T(r,z_i))+O(1)=\ox(T(r))+O(1)
 \]
for all $\a$ such that $p_{i,\a}\ne0$, where $T(r)=\max(T(r,z_i))$.
It follows from
the rules (\ref{Nevanlinna-rules}) that each element $u$ in $\cal M$
satisfies $T(r,u)=\ox(T(r))+O(1)$. Since the rank of $\cal M=n$,
for each $i$ there is an integer 
$d_i$ such that $Z^{d_i}\in\cal M$.
Therefore,
$T(r)=\max T(r,z_i)=\ox(T(r))+O(1)$. This implies $T(r,z_i)=\Ox(1)$
for all $i$, so $z_i$ is constant for each $i$, and by (\ref{ode}),
$p_i(z)=0$ for all $i$. This is impossible since we chose $c=z(0)$
such that (in particular) $p_i(c)\ne0$ for all $i$.

We now take the case $\rank\cal M<n$. Assume first that (\ref{new-ode}) is
complete. Then $w_i$ are 
entire functions and $z_i$ satisfies
 \[
\dot z_i/z_i=p_i(z)=f_i(w)
 \]
The right hand sides are entire functions so integration gives 
 \[
z_i(t)=e^{\int f_i(t)dt},
 \]
so (\ref{ode}) is complete.

If ("ode") is complete, set $w_i=\prod z_j^{a_{ij}}$. These are then
entire function and they satisfy
 \[
\frac{\dot w_i}{w_i}=\sum_j a_{ij}\frac{\dot z_j}{z_j}=\sum_j a_{ij}p_j(z)
=\sum_j a_{ij}f_j(w).
 \]
This is ("new-ode"), which is therefore complete.
\end{proof}

\begin{corollary}\thmfont All complete rational holomorphic vector fields on
$(\complex^{*})^n$ preserve 
the volume form $\bigwedge dz_i/z_i$.
\end{corollary}
\begin{proof} The proof is by induction. We compute
 \[
\frac{d}{dt}\frac{dz_i}{z_i}=\frac{d\dot z_i\,z_i-\dot z_i\,dz_i}{z_i^2}
=\frac{d(z_ip_i(z))z_i-z_ip_i(z)dz_i}{z_i^2}=d(p_i(z)).
\]
This shows in particular that the result holds for $n=1$. Also, we compute
\[
 \frac{d}{dt}\bigwedge_i \frac{dz_i}{z_i}=
\sum_i\frac{dz_1}{z_1}\wedge\dots\wedge\frac{d}{dt}\frac{dz_i}{z_i}\wedge
\dots\wedge\frac{dz_n}{z_n}=\left(\sum_i z_i\p{p_i}{z_i}\right)
\left(\bigwedge_i\frac{dz_i}{z_i}\right).
\]
Since the field is complete, we have Case~2 of the theorem.
We use in particular
that $p_i(z)=f_i(w)$. We have to
prove that 
$\sum z_i{\partial p_i}/{\partial z_i}(z)=0$,
so we compute
\begin{multline}
\sum_i z_i\p{f_i}{z_i}=\sum_i z_i\sum_j\p{f_i}{w_j}\p{w_j}{z_i}=
\sum_i z_i\sum_j\p{f_i}{w_j}w_j\frac{a_{ji}}{z_i}\\
=\sum_{i,j}\p{f_i}{w_j}w_ja_{ji}=\sum_j\p{(af)_j}{w_j}w_j,
\end{multline}
where $(af)_j=\sum a_{ji}f_i$.
By Theorem~"ode-theorem", ("new-ode") is complete and by the induction
hypothesis, the last expression is $0$. The corollary is proved.
\end{proof}

\begin{corollary}\label{dim2}\thmfont All complete rational holomorphic vector
fields on $\cstar2$ are 
of form
 \begin{equation}\label{form2}
   \begin{aligned}
\dot z_1&=z_1(a_2f(z_1^{a_1}z_2^{a_2})+c_1)\\
\dot z_2&=-z_2(a_1f(z_1^{a_1}z_2^{a_2})+c_2),
   \end{aligned}
 \end{equation}
where $f$ is a Laurent polynomial, $a_1$, $a_2$ are integers and $c_1$, $c_2$
are complex numbers. Conversely, all such vector fields are complete.
\end{corollary}

\begin{proof} We use the notation in the theorem. If the field is
complete, we must have $\dim\cal M=1$, so $p_i(Z)=f_i(W)$ for some Laurent
polynomial $f$ and some monomial $W=Z_1^{a_1}Z_2^{a_2}$. The field
("new-ode") becomes
 \[
\dot w=w(a_1f_1(w)+a_2f_2(w)),
 \]
which is complete if and only if $a_1f_1(w)+a_2f_2(w)$ is constant. If
we write $f_1=a_2f+c_1$ and $f_2=-a_1f+c_2$ we get the corollary.
\end{proof}

\begin{remark} In the same way, we can derive (rather complicated)
formulas for the complete vector fields on $\cstar{n}$ for any $n$.
\end{remark}

\begin{remark} The analogs of Theorem~"ode-theorem" and
Corollary~"dim2" for non-rational fields are false. To get a
counterexample, observe that 
 \begin{equation}\label{counterexample}
   \begin{aligned}
\dot z_1&=0\\
\dot z_2&=-z_2z_1
 \end{aligned}
\end{equation}
is complete by Corollary~"dim2". Also, the mapping 
 \begin{align*}
z_1&=\z_1e^{\z_1\z_2}\\
z_2&=\z_2e^{-\z_1\z_2}
 \end{align*}
is a bijection of $\cstar2$ (it is the time~$1$ flow of the field
$\dot\z_1=\z_1^2\z_2$, $\dot\z_2=-\z_1\z_2^2$). If we express
("counterexample") in the new coordinates we get
 \begin{align*}
\dot\z_1&=\z_1^3\z_2e^{\z_1\z_2}\\
\dot\z_2&=-\z_1\z_2(1+\z_1\z_2)e^{\z_1\z_2},
 \end{align*}
which is not of a form corresponding to ("form2").
\end{remark}

\end{document}